\documentclass{elsart}
\usepackage{amsfonts}
\usepackage{amsmath}

\input{tcilatex}
\begin{document}

\begin{frontmatter}

\title{A F\o lner Invariant for Type $II_1$ Factors}

\author{Jon P. Bannon}

\address{Department of Mathematics, Siena College, Loudonville, NY 12211}

\collab{Mohan Ravichandran}

\address{Department of Mathematics and Statistics, The University of New Hampshire, Durham, NH 03824}

\begin{abstract}
In this article we introduce an isomorphism invariant for type $II_1$ factors using the Connes-F\o lner condition.
We compute bounds of this number for free group factors.
\end{abstract}


\end{frontmatter}

\section*{Introduction}

In a series of papers \cite{MN1}-\cite{N2}, Murray and von Neumann
introduced \textquotedblleft rings of operators\textquotedblright , known
nowadays as von Neumann algebras. To them it was clear that what they were
developing was a theory of quantized groups. Many of the examples in their
original paper come from group algebras. Subsequently, concepts and results
in group theory have been a major source of motivation for the development
of operator algebras. Many of the important operator algebra concepts, such
as amenability, property $T$, etc., come directly from properties of various
groups. In this paper, we are concerned with a certain characterization of
amenability for groups due to F\o lner. Our main aim is to introduce an
isomorphism invariant, motivated by F\o lner's characterization, for an
important class of von Neumann algebras called type $II_{1}$ factors.

Von Neumann himself showed that any von Neumann algebra is a direct sum of
\textquotedblleft simple objects\textquotedblright , called factors. These
are weak-operator closed self-adjoint subalgebras of $B(\mathcal{H})$, the
algebra of all bounded operators on a Hilbert space $\mathcal{H}$, whose
centers consist of only scalar multiples of the identity operator. Factors
are called finite if there is a faithful tracial state on them. Those finite
factors which are finite-dimensional as vector spaces are full matrix
algebras $M_{n}(\mathbb{C})$ $(n=1,2,...)$. Those finite factors which are
infinite-dimensional are called factors of type $II_{1}$. In order to
complete the classification of all factors up to $\ast $-isomorphism, it
remains to classify the factors of type $II_{1}$ (cf. \cite{Co2}).

A factor $\mathcal{M}$ is injective if it is the range of a Banach space
projection $\Phi \in B(B(\mathcal{H}))$, for some Hilbert space $\mathcal{H}$%
. There are few computable nontrivial invariants for type $II_{1}$ factors
in general, but the classification of injective factors is complete \cite{Co}%
. It stands to reason that we should try to use tools from the
classification of injective factors to define isomorphism invariants for
general type $II_{1}$ factors. In this paper, we define an invariant $F$\o $%
l(\mathcal{M})$ that will measure how badly a separable type $II_{1}$ factor 
$\mathcal{M}$ fails to satisfy Connes' F\o lner-type condition (Theorem 5.1
in \cite{Co}). We compute explicit bounds for $F$\o $l(\mathcal{M})$ in the
case where $\mathcal{M}$ is the free group factor $L(\mathbb{F}_{n})$.

The layout of the paper is as follows. In the first section, we present some
background on the F\o lner condition for groups in order to provide some
motivation. In the second section we discuss a F\o lner invariant for
groups. In the third section we give some examples of factors and some
questions that will provide further context. In the fourth and final section
we define the pre-invariant $F$\o $l(\mathcal{M},X)$ for a finite subset $X$
of unitary elements in $\mathcal{M}$, and the invariant $F$\o $l(\mathcal{M}%
) $. We then prove that $F$\o $l(\tbigotimes\limits_{j=1}^{\infty }(L(%
\mathbb{F}_{2}))_{j})>0$, and that for any type $II_{1}$ factor $\mathcal{M}$%
, $F$\o $l(\mathcal{M})\leq 2$. Finally, we prove that $F$\o $l(L(\mathbb{F}%
_{n}),X)\leq \sqrt{2-\frac{2}{n^{2}}}$, where $X=%
\{L_{a_{1}},L_{a_{2}},...,L_{a_{n}}\}$ is the set of standard generators of $%
L(\mathbb{F}_{n})$.

The authors wish to thank Professor Liming Ge for many valuable
conversations and for sharing his insight about type $II_{1}$ factors.

The first author was partially supported by a University of New Hampshire
dissertation fellowship.

\section{F\o lner Conditions}

Let $G$ be a discrete group with identity $e$. Let $\mathbb{C}G$ denote the
complex group algebra of formal linear combinations of elements from $G$.
This is a unital $\ast $-algebra, with involution given by the
conjugate-linear extension of the map $g\mapsto g^{-1}$. A faithful trace
state $\tau _{0}$ is defined on $\mathbb{C}G$ by%
\begin{equation*}
\tau _{0}(\tsum \lambda _{g}g)=\lambda _{e}.
\end{equation*}%
Performing the $GNS$ construction using $\tau _{0}$, we faithfully embed $%
\mathbb{C}G$ as a $\ast $-subalgebra $span\{L_{g}:g\in G\}$ of $B(L^{2}(%
\mathbb{C}G,\tau _{0}))$, where the action of $L_{g}$ on $h\in G\subseteq 
\mathbb{C}G$ is given by left-translation in the group $L_{g}h=gh$. We
define the (left) group von Neumann algebra $L(G)$ as $(\mathbb{C}G)^{\prime
\prime }=\overline{\mathbb{C}G}^{WOT}\subseteq B(L^{2}(\mathbb{C}G,\tau
_{0}))$. If we denote by $\langle \cdot ,\cdot \rangle $ the inner product
in $L^{2}(\mathbb{C}G,\tau _{0})$, then $\tau _{0}(T)=\langle Te,e\rangle $,
and $\tau _{0}$ extends by continuity to a trace $\tau $ on all of $L(G)$.
We also have that $\langle g,h\rangle =\delta _{gh}$, so that $G$ is an
orthonormal basis for $L^{2}(\mathbb{C}G,\tau _{0})$. Clearly, by
identifying the standard orthonormal bases, we see that $L^{2}(\mathbb{C}%
G,\tau _{0})\cong L^{2}(L(G),\tau )\cong l^{2}(G)$, and we may consider $%
L(G)\subseteq B(l^{2}(G))$. If $G$ is an $i.c.c.$ group (the conjugacy class
of every $g\not=e$ in $G$ is an infinite set) then $L(G)$ is a factor of
type $II_{1}$.

A discrete group $G$ is amenable if there exists a state on $l^{\infty }(G)$
which is invariant under the left action of $G$ on $l^{\infty }(G).$ Such a
state will be called an invariant mean on $l^{\infty }(G)$. In \cite{Fol}, F%
\o lner used combinatorial methods to find the following condition on a
countable discrete group $G$, and to prove that this condition holds if and
only if $G$ is amenable: \textit{Given }$\{g_{1},g_{2},...,g_{n}\}\subseteq
G $\textit{\ and }$\varepsilon >0,$\textit{\ there exists a finite,
non-empty set }$U\subseteq G$\textit{\ such that }$\forall j\in
\{1,2,...,n\} $%
\begin{equation*}
\frac{\#((g_{j}U\cup U)\backslash (g_{j}U\cap U))}{\#U}\leq \varepsilon .
\end{equation*}%
In \cite{Nam}, I. Namioka was able to prove, using functional analysis, that
an amenable group satisfies F\o lner's condition. The key ingredient in
Namioka's proof is a theorem of Day (see \cite{Nam}, Theorem 2.2).

The classification of injective factors gives us that any two injective type 
$II_{1}$ factors are $\ast $-isomorphic. Furthermore, there are myriad
invariant properties (see \cite{Gr}) that are equivalent to injectivity of a
type $II_{1}$ factor $\mathcal{M\subseteq }B(\mathcal{H})$ with trace $\tau $%
. One such property is Connes' F\o lner-type condition, found in the
statement of Theorem 5.1 in \cite{Co}: \textit{Given }$%
\{x_{1},x_{2},...,x_{n}\}\subseteq M$\textit{\ and }$\varepsilon >0$\textit{%
, there exists a nonzero finite-rank projection }$e\in B(\mathcal{H})$%
\textit{\ such that }$\forall j\in \{1,2,...,n\}$%
\begin{equation*}
||[x_{j},e]||_{H.S.}\leq \varepsilon ||e||_{H.S.}\text{ and }|\tau (x_{j})-%
\frac{\langle x_{j}e,e\rangle _{H.S.}}{\langle e,e\rangle _{H.S.}}|\leq
\varepsilon .
\end{equation*}%
To elucidate the origin of this condition, note that if $\mathcal{M}$ is
injective then $\rho =\tau \circ \Phi $ defines a state on $B(\mathcal{H})$
with the property $\rho |_{\mathcal{M}}=\tau $. Such a state $\rho $ is
called a hypertrace on $\mathcal{M}$. In the case of the von Neumann algebra 
$L(G)$ of a discrete group $G$ we have that $l^{\infty }(G)$ is embedded in $%
B(l^{2}(G))$ as multiplication operators. We see that if $L(G)$ is injective
then $\tau \circ \Phi |_{l^{\infty }(G)}$ is an invariant mean. Conversely,
given an invariant mean on $l^{\infty }(G)$, an averaging process over $%
R(G)(=L(G)^{\prime })$ can be used to construct a conditional expectation of 
$B(l^{2}(G))$ onto $L(G)$, and hence a hypertrace (see 8.7.24 and 8.7.29 of 
\cite{KR}). This suggests that for a general type $II_{1}$ factor, we may
think of a hypertrace as analogous to an invariant mean. Connes exploited
this analogy to prove that when a type $II_{1}$ factor admits a hypertrace
then the factor satisfies the above F\o lner-type condition. The proof of
this follows Namioka's method of obtaining F\o lner's condition from an
invariant mean on a group.

\section{A F\o lner Invariant for Groups}

In \cite{AVR}, Arzhantseva, Burillo, Lustig, Reeves, Short and Ventura have
defined a group invariant $F$\o $l(G)$ that measures how badly a
finitely-generated discrete group $G$ fails to satisfy the classical F\o %
lner condition. In particular this number satisfies, for a group $G$
generated by $n$ elements, the inequality $0\leq $ $F$\o $l(G)\leq \frac{2n-2%
}{2n-1}$. Also $F$\o $l(G)=0$ whenever $G$ is amenable and $F$\o $l(G)=$ $%
\frac{2n-2}{2n-1}$ if and only if $G=\mathbb{F}_{n}$. The notion of boundary
of a subset of a finitely generated group $G$ generally depends on a given
finite generating subset $X$. Arzhantseva et. al.~define 
\begin{equation*}
F\text{\o }l(G,X)=\inf_{\substack{ A\subseteq G  \\ finite}}\frac{\#\partial
_{X}A}{\#A}
\end{equation*}%
where $\partial _{X}A=\{a\in A|ax\not\in A$ for some$~x\in X^{\pm 1}\}$ is
the interior boundary of $A$ with respect to $X$ in $G$. They go on to
define the universal F\o lner invariant%
\begin{equation*}
F\text{\o }l(G)=\inf_{X}F\text{\o }l(G,X)
\end{equation*}%
where the infimum is taken over all finite generating subsets $X$ of $G$.
They prove that if $F$\o $l(G,X)=0$ for some finite generating set $X$ of $G$%
, then $F$\o $l(G,X^{\prime })=0$ for any other finite generating set $%
X^{\prime }$, and this happens only if $G$ is amenable. Non-amenable
discrete groups for which $F$\o $l(G)=0$ are called \textit{weakly amenable }%
and those for which $F$\o $l(G)\not=0$ are called \textit{uniformly
non-amenable}. In \cite{AVR} it is also proven that groups of both types
exist.

In light of the above results, we define the invariant $F$\o $l(\mathcal{M})$
for a type $II_{1}$ factor with separable predual. We note that the analogy
is not entirely straightforward with the group case. The first major
difference is that we exclusively use unitary elements in the computation of 
$F$\o $l(\mathcal{M})$, to avoid blowing up due to scaling by a constant\ in
the Connes-F\o lner condition. The second major difference is that in a type 
$II_{1}$ factor we can find unitary elements arbitrarily norm-close to the
identity, which implies that the second infimum taken in the group case
would always be zero in the new setting. This, in particular, means that the
invariant we introduce will not provide a satisfactory notion of
weak-amenability for type $II_{1}$ factors.

\section{Some Related Examples of Factors}

For the basics of the theory of operator algebras, we refer the reader to 
\cite{KR}$.$

The first classification result in the theory of type $II_{1}$ factors is
the following, due to Murray and von Neumann\cite{MN3}. It remains one of
the deepest results in the subject.

\begin{thm}
Let $\Pi $ denote the group of those permutations of $\mathbb{Z}$ each of
which permutes only finitely many integers, and let $\mathbb{F}_{n}$ be the
nonabelian free group on $n$ generators. Both of these groups are $i.c.c.$,
and give rise to non-isomorphic type $II_{1}$ factors.
\end{thm}

The number $F$\o $l(\mathcal{M})$ will be zero if and only if the factor $%
\mathcal{M}$ is injective. The main problem is to determine whether or not
the invariant can distinguish between a pair of non-injective type $II_{1}$
factors. We are particularly interested in computing the number in the
following two cases.

\begin{exmp}
Let $B(m,n)=\langle a_{1},...,a_{m}|$ $g^{n}=e\rangle $ denote the free
Burnside group on $m$ generators with exponent $n$. If $m>1$ and $n\geq 665$
is odd, then the centralizer of any nonidentity element in $B(m,n)$ is a
cyclic group of order $n$ (cf. \cite{Ad}). It follows in this case that $%
L(B(m,n))$ is a type $II_{1}$ factor. Also, in \cite{Ad2} it is shown that
if $m>1$ and $n\geq 665$ is odd then $B(m,n)$ is not amenable. It follows
from our earlier discussion that $L(B(m,n))$ cannot be an injective factor.
\end{exmp}

\begin{exmp}
Consider R. Thompson's group $F=\langle x_{0},x_{1},x_{2},...|$ $%
x_{i}^{-1}x_{n}x_{i}=x_{n+1},$ $0\leq i\leq n$ $\forall n\in \mathbb{N}%
\rangle $. It is proven in \cite{Jol} that $F$ is $i.c.c.$, and hence that $%
L(F)$ is a type $II_{1}$ factor. A famous conjecture of Geohegan in 1979
asks if $F$ is a non-amenable group which contains no non-abelian free
subgroup. It was proven by Brin and Squire in 1985 that $F$ contains no
non-abelian free subgroup, but it is still unknown whether or not $F$ is an
amenable group(cf. \cite{Jol}).
\end{exmp}

It should be noted that distinguishing the $\ast $-isomorphism classes of
the above type $II_{1}$ factors is an open problem. The last example is
interesting, since finding a single finite subset $X\subseteq F$ with
respect to which the pre-invariant $F$\o $l(L(F),X)\not=0$ amounts to
showing that $F$ is not amenable.

\section{Main Results}

\subsection{The F\o lner Invariant}

We first collect some basic facts about the Hilbert-Schmidt class.

Let $\mathcal{H}$ be a separable Hilbert space. For a positive operator $%
T\in B(\mathcal{H})$, let $Tr(T)=\dsum\limits_{i=1}^{\infty }\langle
Te_{i},e_{i}\rangle $, where $\{e_{i}\}_{i=1}^{\infty }$ is any orthonormal
basis for $\mathcal{H}$. The Hilbert-Schmidt norm of an operator $T\in B(%
\mathcal{H})$ is given by 
\begin{equation*}
||T||_{H.S.}=Tr(T^{\ast }T)^{1/2}.
\end{equation*}%
We say that $T\in B(\mathcal{H})$ is in the Hilbert-Schmidt class when $%
||T||_{H.S.}<\infty $. The class of all such operators in $B(\mathcal{H})$
may be regarded as a Hilbert space when equipped with the inner product $%
\langle A,B\rangle _{H.S.}=Tr(B^{\ast }A).$

Let $\mathcal{M}$ be a factor of type $II_{1}$ with trace $\tau $ acting
standardly on $\mathcal{H}$ $(=L^{2}(\mathcal{M},\tau ))$, and let $\mathcal{%
U}(\mathcal{M})$ be the unitary group of $\mathcal{M}$. Suppose throughout
that $\mathcal{M}$ has separable predual. Connes proves in \cite{Co} that $%
\mathcal{M}$ is injective if and only if the following condition holds:

\begin{quote}
Given $\{x_{1},x_{2},...,x_{n}\}\subset \mathcal{U}(\mathcal{M})$ and $%
\varepsilon >0$, there exists a nonzero finite-rank projection $e\in B(%
\mathcal{H})$ such that $\forall j\in \{1,2,...,n\}$%
\begin{equation*}
||[x_{j},e]||_{H.S.}\leq \varepsilon ||e||_{H.S.}\text{ and }|\tau (x_{j})-%
\frac{\langle x_{j}e,e\rangle _{H.S.}}{\langle e,e\rangle _{H.S.}}|\leq
\varepsilon .
\end{equation*}
\end{quote}

We call this the Connes-F\o lner condition.

\begin{definition}
Let $\mathcal{M}$ be a factor of type $II_{1}$,and $X=%
\{x_{1},x_{2},...,x_{n}\}$ be a finite subset of $\mathcal{U}(\mathcal{M})$.
We define the property $Q(X,\varepsilon )$ to be \textquotedblleft there
exists a nonzero finite-rank projection $e\in B(\mathcal{H})$ such that $%
\forall j\in \{1,2,...,n\}$, $||[x_{j},e]||_{H.S.}\leq \varepsilon
||e||_{H.S.}$ and $|\tau (x_{j})-\frac{\langle x_{j}e,e\rangle _{H.S.}}{%
\langle e,e\rangle _{H.S.}}|\leq \varepsilon .$"
\end{definition}

\begin{definition}
Let $\mathcal{M}$ be a factor of type $II_{1}$, and $X$ be a finite subset
of $\mathcal{U}(\mathcal{M})$. Define%
\begin{equation*}
F\text{\o }l(\mathcal{M},X)=\inf \{\varepsilon >0:Q(X,\varepsilon )\}\text{.}
\end{equation*}
\end{definition}

\begin{definition}
Let $\mathcal{M}$ be a factor of type $II_{1}$. We define the universal F\o %
lner invariant $F$\o $l(\mathcal{M})=\sup_{X}F$\o $l(\mathcal{M},X)$, where
the supremum is taken over all finite sets $X\subset \mathcal{U}(\mathcal{M}%
) $.
\end{definition}

\begin{remark}
By Theorem 5.2 in \cite{Co}, $\mathcal{M}$ is injective if and only if $F$%
\o\ $l(\mathcal{M})=0$.
\end{remark}

We include the following basic observation about monotonicity.

\begin{proposition}
Let $\mathcal{M}$ be a factor of type $II_{1}$. If $X_{1}$ and $X_{2}$ are
finite subsets of $\mathcal{U}(\mathcal{M})$ that generate $\mathcal{M}$ as
a von Neumann algebra, and $X_{1}\subseteq X_{2}$, then $F$\o $l(\mathcal{M}%
,X_{1})\leq F$\o $l(\mathcal{M},X_{2})$.
\end{proposition}

\begin{proof}
We have that for any $\varepsilon >0$ that $Q(X_{2},\varepsilon )\Rightarrow
Q(X_{1},\varepsilon )$, hence 
\begin{equation*}
\inf \{\varepsilon >0:Q(X_{1},\varepsilon )\}\leq \inf \{\varepsilon
>0:Q(X_{2},\varepsilon )\}.
\end{equation*}
\end{proof}

\subsection{Lower Bounds}

\subsubsection{Positivity of $F$\o $l(\tbigotimes\limits_{j=1}^{\infty }(L(%
\mathbb{F}_{2}))_{j})$}

We review the construction of the ultraproduct of finite factors (cf. \cite%
{Co}).

Let $\mathcal{M}^{(n)}$ be finite factors with traces $\tau _{n}$, and let $%
\tprod \mathcal{M}^{(n)}$ denote their $C^{\ast }$-product, i.e. the $%
C^{\ast }$-algebra of uniformly norm-bounded sequences equipped with
coordinatewise operations and the supremum norm. Viewing the Stone-\v{C}ech
compactification $\beta \mathbb{N}$ as the maximal ideal space of $l^{\infty
}(\mathbb{N},\mathbb{C})$, for each $\omega \in \beta \mathbb{N}$ there
corresponds a multiplicative linear functional $\rho \in (l^{\infty }(%
\mathbb{N},\mathbb{C}))^{\#}$. Given $f\in l^{\infty }(\mathbb{N},\mathbb{C}%
) $, we define $\lim_{n\rightarrow \omega }f\equiv \rho (f)$. Consider a
free ultrafilter $\omega \in \beta \mathbb{N}\backslash \mathbb{N}$. We have
that%
\begin{equation*}
\mathcal{I}_{\omega }=\{(A_{i})_{i}\in \tprod \mathcal{M}^{(n)}:\lim_{i%
\rightarrow \omega }\tau _{i}(A_{i}^{\ast }A_{i})=0\}
\end{equation*}%
is a closed two-sided ideal in $\tprod \mathcal{M}^{(n)}$, and by a result
of Sakai \cite{Sa}, the quotient $(\tprod \mathcal{M}^{(n)})/\mathcal{I}%
_{\omega }$ is a factor von Neumann algebra algebra $\tprod^{\omega }%
\mathcal{M}^{(n)}$ with a faithful, normal trace $\tau _{\omega }$ defined
by $\tau _{\omega }((A_{i})_{i}+\mathcal{I}_{\omega })=\lim_{i\rightarrow
\omega }\tau _{i}(A_{i})$. The factor $\tprod^{\omega }\mathcal{M}^{(n)}$
will be called an ultraproduct of the $\mathcal{M}^{(n)}$ with respect to
the free ultrafilter $\omega $, or simply an ultraproduct of the $\mathcal{M}%
^{(n)}$. If $\mathcal{M}$ is a finite factor and $\mathcal{M}^{(n)}=\mathcal{%
M}$ for all $n$, then the ultraproduct is called an ultrapower, and is
written as $\mathcal{M}^{\omega }$. In this case, we embed $\mathcal{M}$ in $%
\mathcal{M}^{\omega }$ as constant sequences.

In what follows, let $\tau _{k}$ denote the normalized trace on the
appropriate type $I_{k}$ factor.

\begin{lemma}
Suppose that $\mathcal{M}$ is a type $II_{1}$ factor with trace $\tau $. If $X$ is a finite subset of $\mathcal{U}(\mathcal{M})$ and $F$\o $l(\mathcal{%
M},X)=0$, then for every $U\in X$, 
$M \in \mathbb{N}$ and $\delta>0$ there exists $m\in \mathbb{N}$ such that $m \geq M$,  $Q(X,\frac{1}{m} )$ via a projection $e_{m}$ of some finite rank $l(m)$ and
there exists a unitary element $W_{m}\in e_{m}\mathcal{B}(L^{2}(\mathcal{M},\tau))e_{m}$ satisfying
\begin{equation*}
||e_{m}Ue_{m}-W_m||_{\tau _{l(m)}}\leq \delta \text{.}
\end{equation*} 
\end{lemma}

\begin{proof}
Let $M>0$. If $F$\o $l(\mathcal{M},X)=0$, then there is an increasing sequence of positive integers $(n_{k})_{k \in \mathbb{N}}$
such that $Q(X,\frac{1}{n_{k}}),$ and hence there is a projection $e_{n_{k}}$ of some finite rank $l(n_{k})$ onto the span of an orthonormal set $\{\xi
_{i}^{(n_k)}\}_{i=1}^{l(n_k)}$ of vectors in $L^{2}(\mathcal{M})$ satisfying%
\begin{equation*}
0\leq \frac{||[U,e_{n_{k}}]||_{H.S.}}{||e_{n_{k}}||_{H.S.}}=\sqrt{2}\sqrt{%
1-||e_{n_{k}}Ue_{n_{k}}||_{\tau _{l(n_{k})}}^{2}}\leq \frac{1}{n_{k}}
\end{equation*}%
for all $U\in X$. With $e_{n_{k}}Ue_{n_{k}}=A_{n_{k}}=[\langle \xi _{q}^{(n_{k})},U\xi
_{p}^{(n_{k})}\rangle ]_{q,p=1}^{l(n_{k})}$, we have 
\begin{equation*}
1-\frac{1}{2n_{k}^{2}}\leq \tau _{l(n_{k})}(A_{n_{k}}^{\ast }A_{n_{k}})=\tau
_{l(n_{k})}(A_{n_{k}}A_{n_{k}}^{\ast })=||e_{n_{k}}Ue_{n_{k}}||_{\tau _{l(n_{k})}}^{2}.
\end{equation*}

Furthermore, since $e_{n_{k}}$ is a projection, $||e_{n_{k}}||\leq 1$ and hence%
\begin{equation*}
||A_{n_{k}}||=||e_{n_{k}}Ue_{n_{k}}||\leq ||U||||e_{n_{k}}||^{2}\leq 1,
\end{equation*}%
and hence $||A_{n_{k}}A_{n_{k}}^{\ast }||=||A_{n_{k}}||^{2}\leq 1$. Let $\omega $ be a
free ultrafilter, and $\tprod\limits^{\omega }M_{l(n_{k})}(\mathbb{C})$ denote
the ultraproduct factor as defined above. We have a sequence $(A_{n_{k}})=(A_{n_{k}})_{n_{k}\geq M}$
of matrices satisfying 
\begin{eqnarray*}
\tau _{\omega }((A_{n_{k}}^{\ast }A_{n_{k}})+\mathcal{I}_{\omega }) &=&\tau _{\omega
}((A_{n_{k}}A_{n_{k}}^{\ast })+\mathcal{I}_{\omega })=1 \\
&=&\tau _{\omega }((I_{n_{k}})+\mathcal{I}_{\omega })
\end{eqnarray*}%
so by faithfulness of $\tau _{\omega }$ and the fact that $%
(I_{n_{k}}-A_{n_{k}}A_{n_{k}}^{\ast })_{n_{k}}\geq 0$ for all $n$,%
\begin{equation*}
\tau _{\omega }((I_{n_{k}}-A_{n_{k}}A_{n_{k}}^{\ast })_{n_{k}}+\mathcal{I}_{\omega })=0\text{%
,}
\end{equation*}%
so indeed $(A_{n_{k}})$ represents a unitary element in the ultraproduct $%
\tprod\limits^{\omega }M_{l(n_{k})}(\mathbb{C})$. Recall that if 
\begin{equation*}
(A_{n_{k}})+\mathcal{I}_{\omega }\text{ and }(B_{n_{k}})+\mathcal{I}_{\omega }
\end{equation*}%
represent distinct elements of $\tprod\limits^{\omega }M_{l(n_{k})}(\mathbb{C})$%
, then the $2$-norm distance between them is given by%
\begin{eqnarray*}
&&||(A_{n_{k}}-B_{n_{k}})+\mathcal{I}_{\omega }||_{2} \\
&=&\tau _{\omega }(((A_{n_{k}}^{\ast }-B_{n_{k}}^{\ast })+\mathcal{I}_{\omega
})((A_{n_{k}}-B_{n_{k}})+\mathcal{I}_{\omega }))^{1/2} \\
&=&[\lim_{l(n_{k})\rightarrow \omega }\tau _{l(n_{k})}((A_{n_{k}}^{\ast }-B_{n_{k}}^{\ast
})(A_{n_{k}}-B_{n_{k}}))]^{1/2} \\
&=&[\lim_{l(n_{k})\rightarrow \omega }||A_{n_{k}}-B_{n_{k}}||_{\tau _{l(n_{k})}}^{2}]^{1/2}.
\end{eqnarray*}%
Suppose that $\delta >0$ and that for every unitary $l(n_{k})\times l(n_{k})$ matrix 
$W_{n_{k}}$, $||A_{n_{k}}-W_{n_{k}}||_{\tau _{l(n_{k})}}>\delta ,$ it then follows that $%
||(A_{n_{k}}-W_{n_{k}})+\mathcal{I}_{\omega }||_{2}>\delta $ in $L^{2}(\tprod%
\limits^{\omega }M_{l(n_{k})}(\mathbb{C}),\tau _{\omega })$. Since every
sequence $(W_{n_{k}})$ represents a unitary element in $\tprod\limits^{\omega
}M_{l(n_{k})}(\mathbb{C})$, and from the polar decomposition and the fact that the ultraproduct is a finite factor every 
unitary element is represented by such a
sequence, a contradiction follows, since $(A_{n_{k}})$ represents a unitary
element in $\tprod\limits^{\omega }M_{l(n_{k})}(\mathbb{C})$. Therefore, for all 
$\delta >0$ there exists a unitary $l(n_{k})\times l(n_{k})$ matrix $W_{n_{k}}$ so that $%
||A_{n_{k}}-W_{n_{k}}||_{\tau _{l(n_{k})}}\leq \delta $, hence we may view $W_{n_{k}}$ as a
unitary element of $e_{n_{k}}B(L^{2}(\mathcal{M}))e_{n_{k}}$ (i.e. a unitary
operator on $span\{\xi _{i}^{(n_{k})}\}_{i=1}^{l(n_{k})}\cong \mathbb{C}^{l(n_{k})}$).
\end{proof}

We recall the construction of the infinite tensor product of a collection of
finite factors. Let $\{\mathcal{M}_{i}\}_{i\in \mathbb{N}}$ be a countable
collection of finite factors with faithful normal traces $\tau _{i}$, and
let $\mathcal{A}_{n}\equiv \tbigotimes\limits_{i=1}^{n}\mathcal{M}_{i}$
denote an algebraic tensor product. The map $T_{1}\otimes ...\otimes
T_{n}\mapsto T_{1}\otimes ...\otimes T_{n}\otimes I$ on simple tensors
extends to a unital embedding of $\mathcal{A}_{n}$ into $\mathcal{A}_{n+1}$.
Let $\mathcal{A}$ be the direct limit algebra obtained via these embeddings.
We have that $\mathcal{A}$ obtains a unital $\ast $-algebra structure and a
faithful normal trace $\tau _{0}$ from the $\mathcal{M}_{i}$. Let $\pi $
denote the GNS representation obtained from $\mathcal{A}$ and $\tau _{0}$.
We define $\tbigotimes\limits_{i=1}^{\infty }\mathcal{M}_{i}\equiv \pi (%
\mathcal{A})^{\prime \prime }$. It is easy to see that this is a factor. The
state $\tau _{0}$ extends uniquely to a faithful normal trace on $%
\tbigotimes\limits_{i=1}^{\infty }\mathcal{M}_{i}$, so we obtain that the
factor is finite.

The central sequence algebra $\mathcal{M}_{\omega }=\mathcal{M}^{\prime
}\cap \mathcal{M}^{\omega }$ is the algebra of all elements in $\mathcal{M}%
^{\omega }$ that commute with $\mathcal{M}$ (see \cite{Dix}, \cite{McD}, 
\cite{Co}). If $\mathcal{M}_{\omega }\not=\mathbb{C}I$, then we say that $%
\mathcal{M}$ has property $\Gamma $. It is a straightforward exercise to
show that every infinite tensor product factor $\tbigotimes\limits_{i=1}^{%
\infty }\mathcal{M}_{i}$ has property $\Gamma $.

In the next theorem, let $\mathcal{M}$ denote the type $II_{1}$ factor $%
\tbigotimes\limits_{j=1}^{\infty }(L(\mathbb{F}_{2}))_{j}$, and let $%
U=L_{a}\otimes I\otimes I...$ and $V=L_{b}\otimes I\otimes I...$ in $%
\mathcal{M}$. We now compute an explicit lower bound for $F$\o $l(\mathcal{M}%
)$.

\begin{theorem}
If $X=\{U,V\}$, then $F$\o $l(\mathcal{M},X)>0$.
\end{theorem}

\begin{proof}
Suppose that $F$\o $l(\mathcal{M},X)=0$, so by the lemma, there exists a positive integer $n$ and
a rank $n$ projection $e\in B(L^{2}(\mathcal{M}))$ such that%
\begin{equation*}
0\leq \frac{||[U,e]||_{H.S.}}{||e||_{H.S.}}=\sqrt{2}\sqrt{1-||eUe||_{\tau
_{n}}^{2}}\leq \frac{1}{7}.
\end{equation*}%
We have that $||eUe-Ue||_{\tau _{n}}^{2}=1-||eUe||_{\tau _{n}}^{2}\leq \frac{%
1}{98}.$ By the above lemma, there is an $n\times n$ unitary matrix $W\in
eB(L^{2}(\mathcal{M}))e$ such that 
\begin{equation*}
||eUe-W||_{\tau _{n}}\leq (1-\frac{1}{\sqrt{2}})\frac{1}{7}.
\end{equation*}

By the triangle inequality, we have that%
\begin{equation*}
||Ue-W||_{\tau _{n}}\leq \frac{1}{7}.
\end{equation*}%
Let $\{\xi _{1},...,\xi _{n}\}\subseteq L^{2}(\mathcal{M})$ be an
orthonormal basis for the range of $e$. Since $W\in eB(L^{2}(\mathcal{M}))e$
we have%
\begin{equation*}
||Ue-W||_{\tau _{n}}^{2}=\frac{1}{n}\sum_{i=1}^{n}||(Ue-W)\xi _{i}||^{2}.
\end{equation*}%
Writing $(g_{1},g_{2},...)\in \mathbb{F}_{2}\times \mathbb{F}_{2}...\times 
\mathbb{F}_{2}...$ in place of $\chi _{\{g_{1}\}}\otimes $ $\chi
_{\{g_{2}\}}\otimes ...$, we may view $\mathbb{F}_{2}^{\infty }\mathbb{=F}%
_{2}\times \mathbb{F}_{2}...\times \mathbb{F}_{2}...$ as an orthonormal
basis for $L^{2}(\mathbb{F}_{2}^{\infty })\cong L^{2}(\mathcal{M})$. Consider 
the action of $\mathbb{F}_{2}$ on $\mathbb{F}_{2}^{\infty }$ in the first
coordinate, that is, the action $g\in \mathbb{F}_{2}$ given by $%
g(g_{1},g_{2},...)=(gg_{1},g_{2},...)$. For $i\in \{1,2,...,n\}$, if 
\begin{equation*}
\xi_{i}=\sum_{(g_{1},g_{2},...)\in \mathbb{F}_{2}^{\infty }}\lambda
_{(g_{1},g_{2},...)}^{(i)}(g_{1},g_{2},...)
\end{equation*} %
then%
\begin{eqnarray*}
||(Ue-W)\xi _{i}||^{2} &=&||(U-W)\sum_{(g_{1},g_{2},...,)\in \mathbb{F}%
_{2}^{\infty }}\lambda _{(g_{1},g_{2},...)}^{(i)}(g_{1},g_{2},...)||^{2} \\
&=&\sum_{(g_{1},g_{2},...)\in \mathbb{F}_{2}^{\infty }}|\lambda
_{(a^{-1}g_{1},g_{2},...)}^{(i)}-\sum_{k=1}^{n}W_{ik}\lambda
_{(g_{1},g_{2},...)}^{(k)}|^{2}\text{.}
\end{eqnarray*}%
For $S$ a non-empty subset of $\mathbb{F}_{2}^{\infty }$ and 
\begin{equation*}
\eta =\sum_{(g_{1},g_{2},...)\in \mathbb{F}_{2}^{\infty }}\mu
_{(g_{1},g_{2},...)}(g_{1},g_{2},...)\in L^{2}(\mathbb{F}_{2}^{\infty }),
\end{equation*}%
define 
\begin{equation*}
||\eta ||_{S}^{2}\equiv \sum_{(g_{1},g_{2},...)\in S}|\mu
_{(g_{1},g_{2},...)}|^{2}.
\end{equation*}%
It follows that%
\begin{equation*}
||(U-W)\xi _{i}||_{S}^{2}=\sum_{(g_{1},g_{2},...)\in S}|\lambda
_{(a^{-1}g_{1},g_{2},...)}^{(i)}-\sum_{k=1}^{n}W_{ik}\lambda
_{(g_{1},g_{2},...)}^{(k)}|^{2}.
\end{equation*}%
We have that%
\begin{eqnarray*}
|(||U\xi _{i}||_{S}-||W\xi _{i}||_{S})| &\leq &||(U-W)\xi _{i}||_{S} \\
&\leq &||(U-W)\xi _{i}||.
\end{eqnarray*}%
and using the inequality $(x_{1}+...+x_{n})^{2}\leq
n(x_{1}^{2}+...+x_{n}^{2})$ and the triangle inequality, we get%
\begin{eqnarray*}
&&|\frac{1}{n}\sum_{i=1}^{n}||U\xi _{i}||_{S}^{2}-\frac{1}{n}%
\sum_{i=1}^{n}||W\xi _{i}||_{S}^{2})|^{2} \\
&\leq &\frac{n}{n^{2}}\sum_{i=1}^{n}|(||U\xi _{i}||_{S}^{2}-||W\xi
_{i}||_{S}^{2})|^{2} \\
&\leq &\frac{1}{n}\sum_{i=1}^{n}||(||U\xi _{i}||_{S}-||W\xi
_{i}||_{S})|(||U\xi _{i}||_{S}+||W\xi _{i}||_{S})|^{2} \\
&\leq &\frac{4}{n}\sum_{i=1}^{n}|(||U\xi _{i}||_{S}-||W\xi _{i}||_{S})|^{2}
\\
&\leq &\frac{4}{n}\sum_{i=1}^{n}||(U-W)\xi _{i}||_{S}^{2} \\
&\leq &\frac{4}{n}\sum_{i=1}^{n}||(U-W)\xi _{i}||^{2} \\
&\leq &\frac{4}{49}.
\end{eqnarray*}%
With $\eta =\sum_{(g_{1},g_{2},...)\in \mathbb{F}_{2}^{\infty }}\mu
_{(g_{1},g_{2},...)}(g_{1},g_{2},...)\in L^{2}(\mathbb{F}_{2}^{\infty })$,
define 
\begin{equation*}
\eta |_{S}\equiv \sum_{(g_{1},g_{2},...)\in S}\mu
_{(g_{1},g_{2},...)}(g_{1},g_{2},...)\in L^{2}(S)\subseteq L^{2}(\mathbb{F}%
_{2}^{\infty }).
\end{equation*}%
Note that $||\eta |_{S}||_{L^{2}(S)}=||\eta ||_{S}$. We have that%
\begin{eqnarray*}
W\xi _{i}|_{S} &=&\sum_{k=1}^{n}W_{ik}\xi _{k}|_{S}=\sum_{g\in
S}(\sum_{k=1}^{n}W_{ik}\lambda _{g}^{(k)})g \\
&=&(W\xi _{i})|_{S}\text{,}
\end{eqnarray*}%
We may conclude, since $W$ is a unitary operator on $\mathbb{C}^{n}$, that%
\begin{eqnarray*}
\sum_{i=1}^{n}||W\xi _{i}||_{S}^{2} &=&\sum_{i=1}^{n}||(W\xi
_{i})|_{S}||_{L^{2}(S)}^{2} \\
&=&\sum_{i=1}^{n}||W\xi _{i}|_{S}||_{L^{2}(S)}^{2} \\
&=&\sum_{i=1}^{n}||\xi _{i}|_{S}||_{L^{2}(S)}^{2}=\sum_{i=1}^{n}||\xi
_{i}||_{S}^{2}\text{.}
\end{eqnarray*}%
We also have that for each $i$,%
\begin{eqnarray*}
(U\xi _{i})|_{S} &=&\sum_{(g_{1},g_{2},...)\in S}\lambda
_{(a^{-1}g_{1},g_{2},...)}^{(i)}(g_{1},g_{2},...) \\
&=&U\xi _{i}|_{S}\text{.}
\end{eqnarray*}

It follows that%
\begin{eqnarray*}
\sum_{i=1}^{n}||U\xi _{i}|_{S}||_{L^{2}(S)}^{2} &=&\sum_{i=1}^{n}||(U\xi
_{i})|_{S}||_{L^{2}(S)}^{2} \\
&=&\sum_{i=1}^{n}||U\xi _{i}||_{S}^{2}.
\end{eqnarray*}%
Notice that%
\begin{eqnarray*}
||U\xi _{i}||_{S}^{2} &=&||(U\xi _{i})|_{S}||_{L^{2}(S)}^{2} \\
&=&\sum_{(g_{1},g_{2},...)\in S}|\lambda
_{(a^{-1}g_{1},g_{2},...)}^{(i)}|^{2} \\
&=&\sum_{(g_{1},g_{2},...)\in a^{-1}S}|\lambda
_{(g_{1},g_{2},...)}^{(i)}|^{2} \\
&=&||\xi _{i}|_{a^{-1}S}||_{L^{2}(a^{-1}S)}^{2}=||\xi _{i}||_{a^{-1}S}^{2}.
\end{eqnarray*}%
We have that 
\begin{eqnarray*}
|\frac{1}{n}\sum_{i=1}^{n}(||\xi _{i}||_{a^{-1}S}^{2}-||\xi _{i}||_{S}^{2})|
&=&|\frac{1}{n}\sum_{i=1}^{n}||U\xi _{i}||_{S}^{2}-\frac{1}{n}%
\sum_{i=1}^{n}||\xi _{i}||_{S}^{2})|^{2} \\
&\leq &\frac{4}{49}.
\end{eqnarray*}

Now we shall choose a subset $S$ for which the above inequality will give us
a contradiction. For simplicity of notation, let us define%
\begin{equation*}
c_{S}\equiv \frac{1}{n}\sum_{i=1}^{n}||\xi _{i}||_{S}^{2}\text{.}
\end{equation*}%
The above inequality becomes%
\begin{equation*}
|c_{a^{-1}S}-c_{S}|\leq \frac{4}{49}\text{.}
\end{equation*}%
If we carry out the above analysis using $V$ in place of $U$, we obtain%
\begin{equation*}
|c_{b^{-1}S}-c_{S}|\leq \frac{4}{49}.
\end{equation*}%
Since $S$ was arbitrary, we could replace $S$ by $aS$ (resp. $bS$) to get%
\begin{equation*}
|c_{S}-c_{aS}|\leq \frac{4}{49}
\end{equation*}%
\begin{equation*}
(\text{resp. }|c_{S}-c_{bS}|\leq \frac{4}{49}).
\end{equation*}%
Choose the set $S$ to be $S^{\prime }\times $ $\mathbb{F}_{2}^{\infty }$,
where $S^{\prime }$ is the set of all reduced words in $\mathbb{F}_{2}$ that
begin with $a^{-1}$. Then $S\cup aS=\mathbb{F}_{2}^{\infty }$ and also $S,$ $%
bS$ and $b^{-1}S$ are pairwise disjoint. Since $S\cup aS=\mathbb{F}%
_{2}^{\infty }$, we have that $c_{S}$ or $c_{aS}$ exceeds $\frac{1}{2}$.
Since $S,$ $bS$ and $b^{-1}S$ are pairwise disjoint, at least one of $%
c_{S},c_{bS}$ or $c_{b^{-1}S}$ must be $\leq \frac{1}{3}$. With no loss of
generality, we may assume that $\frac{1}{2}\leq c_{aS}$. It follows that 
\begin{equation*}
\frac{1}{2}\leq c_{aS}\leq |c_{S}-c_{aS}|+|c_{S}|\leq \frac{4}{49}+c_{S}
\end{equation*}%
so that%
\begin{equation*}
\frac{1}{2}-\frac{4}{49}\leq c_{S}\text{.}
\end{equation*}%
Let us assume, again with no loss of generality, that $c_{bS}\leq \frac{1}{3}
$, then 
\begin{equation*}
c_{S}\leq |c_{S}-c_{bS}|+c_{bS}\leq \frac{4}{49}+\frac{1}{3}\text{.}
\end{equation*}%
It follows that%
\begin{equation*}
\frac{5}{12}<\frac{1}{2}-\frac{4}{49}\leq c_{S}\leq \frac{1}{3}+\frac{4}{49}<%
\frac{5}{12}
\end{equation*}%
which is a contradiction.
\end{proof}

\begin{remark}
The above proof, slightly modified, gives that
$$F\text{\o } l(L(\mathbb{F}_{2}),\{L_{a},L_{b}\})>0.$$
\end{remark}

\subsection{Upper Bounds}

We begin this section by proving that the universal F\o lner constant of any
given type $II_{1}$ factor cannot exceed $2.$ We then move on to compute
specific upper bounds for $F$\o $l(L(\mathbb{F}_{n}),X)$, with $X$ the set
of standard generators.

\begin{proposition}
For any type $II_{1}$ factor $\mathcal{M}$, $F$\o $l(\mathcal{M})\leq 2$.
\end{proposition}

\begin{proof}
First suppose that $X$ is a finite set of unitary elements in $M$, such that 
$\varepsilon >2$ and the negation of $Q(X,\varepsilon )$ holds. If $k\mathbb{\in N%
}$ and $e$ is a rank $k$ projection such that $\sqrt{2}\sqrt{1-||eUe||_{\tau
_{k}}^{2}}>\varepsilon $ then $||eUe||_{\tau _{k}}^{2}<0$, which cannot
happen. It follows that for every $k\mathbb{\in N}$ and rank $k$ projection $%
e$ in $B(L^{2}(M))$, there exists $U\in X$ such that 
\begin{equation*}
|\tau (U)-\tau _{k}(eUe)|>\varepsilon \text{.}
\end{equation*}%
However, using the triangle and Cauchy-Schwartz inequalities, 
\begin{equation*}
2<\varepsilon <|\tau (U)-\tau _{k}(eUe)|\leq |\tau (U)|+|\tau _{k}(eUe)|\leq
2\text{, }
\end{equation*}%
a contradiction.
\end{proof}

\begin{proposition}
$F$\o $l(L(\mathbb{F}_{n}),X)\leq \sqrt{2-\frac{2}{n^{2}}}$, where 
\begin{equation*}
X=\{L_{a_{1}},L_{a_{2}},...,L_{a_{n}}\}
\end{equation*}
 is the set of standard generators
of $L(\mathbb{F}_{n})$.
\end{proposition}

\begin{proof}
For $i\in \{1,2,...,n\}$ and $\varepsilon \in \{\pm 1\}$, define 
\begin{equation*}
S_{a_{i}^{\varepsilon }}=\{g\in \mathbb{F}_{n}|\text{ }g\text{ begins with }%
a_{i}^{\varepsilon }\}.
\end{equation*}%
Given $i\in \{1,2,...,n\}$, Let $w$ be the least positive integer equivalent
to $(i-1)$ modulo $n,$ and $\{g_{j}^{(i)}|$ $j\in \mathbb{N\}}$ be the list
of elements in $S_{a_{w}^{-1}}$. For $m\in \{1,2,...,k\}$, let 
\begin{equation*}
\xi _{m}=\sum_{t=1}^{\infty }\frac{1}{\sqrt{(n+1)^{t}}}%
(\sum_{i=1}^{n}a_{i}^{m}g_{t}^{(i)})\in L^{2}(\mathbb{F}_{n}).
\end{equation*}%
We have $\{\xi _{m}\}_{m=1}^{k}$ is an orthonormal set.

Let $e$ be the projection onto $span\{\xi _{m}\}_{m=1}^{k}$. We have that
for all $j\in \{1,2,...,n\}$ and $m,s\in \{1,2,...,k\}$ that 
\begin{equation*}
\langle L_{a_{j}}\xi _{s},\xi _{m}\rangle =0
\end{equation*}%
unless $m>1$ and $s=m-1$, in this case%
\begin{equation*}
\langle L_{a_{j}}\xi _{m-1},\xi _{m}\rangle =\frac{1}{n}\text{.}
\end{equation*}%
It follows that $||eL_{a_{j}}e||_{\tau _{k}}^{2}=\frac{1}{k}%
\sum_{m=2}^{k}|\langle L_{a_{j}}\xi _{m-1},\xi _{m}\rangle |^{2}=\frac{k-1}{%
kn^{2}}$, and hence 
\begin{equation*}
\sqrt{2}\sqrt{1-||eL_{a_{j}}e||_{\tau _{k}}^{2}}=\sqrt{2}\sqrt{1-(\frac{k-1}{%
k})\frac{1}{n^{2}}}.
\end{equation*}%
It follows that $F$\o $l(L(\mathbb{F}_{n}),X)\leq \inf_{k\in \mathbb{N}}\{%
\sqrt{2}\sqrt{1-(\frac{k-1}{k})\frac{1}{n^{2}}}\}=\sqrt{2-\frac{2}{n^{2}}}$.
\end{proof}

\end{document}